\documentclass[12pt]{amsart}
\usepackage{amssymb}

\newtheorem{theorem}{Theorem}[section]
\newtheorem{thm}[theorem]{Theorem}
\newtheorem{lemma}[theorem]{Lemma}
\newtheorem{obs}[theorem]{Observation}
\newtheorem{cor}[theorem]{Corollary}
\newtheorem*{conj}{Conjecture}
\theoremstyle{definition}
\newtheorem{defn}[theorem]{Definition}

\def\s{\subseteq}
\def\bks{\setminus}
\def\gorer{\Longrightarrow}
\def\pd{\Game}
\DeclareMathOperator{\cf}{cf}
\DeclareMathOperator{\im}{Im}
\DeclareMathOperator{\cov}{cov}

\begin{document}
\title[Antichains in posets of singular cofinality]{Antichains in partially ordered sets of singular cofinality}

\author{Assaf Rinot}
\address{School of Mathematical Sciences\\Tel Aviv University\\Tel Aviv 69978, Israel}
\email{paper03@rinot.com}
\urladdr{http://www.tau.ac.il/~rinot}


\begin{abstract}In their paper from 1981, Milner and Sauer conjectured that for any poset $\langle P,\le\rangle$, if $\cf(P,\le)=\lambda>\cf(\lambda)=\kappa$,
then $P$ must contain an antichain of size $\kappa$.

We prove that for $\lambda>\cf(\lambda)=\kappa$, if there exists a cardinal $\mu<\lambda$ such that $\cov(\lambda,\mu,\kappa,2)=\lambda$,
then any poset of cofinality $\lambda$ contains $\lambda^\kappa$ antichains of size $\kappa$.

The hypothesis of our theorem is very weak and is a consequence of many well-known axioms such as GCH, SSH and PFA.
The consistency of the negation of this hypothesis is unknown.
\end{abstract}

\maketitle
\section{Introduction}
\subsection{Background}
Assume $\langle P,\le \rangle$ is a poset.
For $A\s P$, let the \emph{downward closure} of $A$ be $\underline{A}:=\{ x\in P\mid \exists y\in A(x\le y)\}$,
the \emph{upward closure} of $A$ be $\overline{A}:=\{ x\in P\mid \exists y\in A(y\le x)\}$,
the \emph{external cofinality of $A$} be  $\cf_{P}(A):=\min\{ |B| \mid B\subseteq P, A\subseteq \underline{B}\}$,
and the \emph{cofinality} of the whole poset be $\cf(P,\le)=\cf_P(P)$.
If $P\s\underline{A}$,  we say that $A$ is \emph{cofinal} in $P$.

For $x,y\in P$, we say that $x$ and $y$ are \emph{incomparable} iff $x\not\le y$ and $y\not\le x$.
$A\s P$ is said to be an \emph{antichain} iff $x,y$ are {\bf incomparable} for all distinct $x,y\in A$.

In his paper \cite{pouzet}, Pouzet proved his celebrated 
theorem
stating that any updirected poset with no infinite antichain contains
a cofinal subset which is isomorphic to a product of finitely many regular cardinals.

Since any poset with no infinite antichain is the union of finitely many updirected subposets,
we have:
\begin{thm}[Pouzet \cite{pouzet}] Assume $\langle P,\le\rangle$ is a poset.

If $\cf(P,\le)$ is a singular cardinal, then $P$ contains an infinite antichain.
\end{thm}
This lead to the formulation of a very natural conjecture,
first appearing implicitly in \cite{pouzet}, and then explicitly in \cite{MiSa}:
\begin{conj}[Milner-Sauer \cite{MiSa}]Assume $\langle P,\le\rangle$ is a poset.

If $\cf(P,\le)=\lambda>\cf(\lambda)=\kappa$, then $P$ contains an antichain of size $\kappa$.
\end{conj}
This conjecture and further generalizations of it were the subject of research of
\cite{gorelic,MR1197093,HaSa,magidor,MR779842,MR661297,MR835967,pouzet2,MR699097,rinot2,rinot}.
For $\lambda>\cf(\lambda)=\kappa$, Milner and Prikry \cite{MR699097} proved that $\mu^{<\kappa}<\lambda$ for all $\mu<\lambda$,
implies that any poset of cofinality
$\lambda$ indeed contains an antichain of size $\kappa$.
Milner and Pouzet \cite{MR661297} derived the same result already from $\lambda^{<\kappa}=\lambda$.
Hajnal and Sauer \cite{HaSa} obtained $\lambda^\kappa$ antichains (of size $\kappa$),
whenever $\lambda$ is a (singular) strong limit,
and this was later improved in Milner and Pouzet \cite{pouzet2}, and Gorelic \cite{gorelic},
yielding $\lambda^\kappa$ antichains already from  $\lambda^{<\kappa}=\lambda$.

The current state of the conjecture is the following:
\begin{thm}[\cite{rinot}]\label{11} Assume cardinals $\lambda>\cf(\lambda)=\kappa$.

If $\cf([\lambda]^{<\kappa},\s)=\lambda$, then any poset of cofinality $\lambda$ contains $\lambda^\kappa$ antichains of size $\kappa$.\footnote{We
were informed by U. Abraham that  M. Magidor \cite{magidor} independently obtained a short combinatorial proof
that if $\cf([\lambda]^{<\kappa},\s)=\lambda$, then any poset of cofinality $\lambda$ contains an antichain of size $\kappa$.}
\end{thm}
The main difference between the hypothesis $\lambda^{<\kappa}=\lambda$ and $\cf([\lambda]^{<\kappa},\s)=\lambda$
is that the first can easily be violated using, e.g., Cohen forcing,
while large cardinals are necessary for the violation of the second hypothesis (Cf. \cite{rinot}).

In this paper, we improve Theorem \ref{11} to the following:

\begin{thm}\label{13} Assume cardinals $\lambda>\cf(\lambda)=\kappa$.

If there exists a cardinal $\mu<\lambda$ such that $\cov(\lambda,\mu,\kappa,2)=\lambda$,
then any poset of cofinality $\lambda$ contains $\lambda^\kappa$ antichains of size $\kappa$.\footnote{For the notion of $\cov(\lambda,\mu,\kappa,\sigma)$, see Definition \ref{covdef}.}
\end{thm}
To appreciate the improvement, we mention that while the negation of the hypothesis of Theorem \ref{11}
can indeed be obtained via forcing with large cardinals, the consistency of the negation
of the latter hypothesis is unknown.
Presenting a model with $\cov(\lambda,\mu,\kappa,2)>\lambda$ for all $\mu\in(\kappa,\lambda)$
is one of the basic open problems of modern cardinal arithmetic.

It is also worth mentioning that a crucial part in the proof of Theorem \ref{11} in \cite{rinot}
was metamathematical, that is,
Gitik's theorem that $\cf([\lambda]^{<\kappa},\s)=\lambda$ implies $L[A]\models\lambda^{<\kappa}=\lambda$
for a particular relevant subset $A\s\lambda$.

In this paper, by extending the methods of \cite{gorelic,HaSa,pouzet2}, a purely combinatorial proof is obtained.%
\footnote{But to the topological result of \cite{rinot} about spaces of singular density, whose proof uses Gitik's theorem,
there is no purely combinatorial proof that we know of.}

\subsection{Notation}
We denote cardinals with the Greek letters, $\lambda,\kappa,\mu,\theta,\sigma$
and ordinals with the letters $\alpha,\beta,\gamma,\delta,\tau$.
For a set $A$, a cardinal $\mu$, and a binary relation $\lhd\in\{<,\le\}$, let $[A]^{\lhd\mu}:=\{ X\s A \mid |A|\lhd\mu\}$,
and $[A]^{\mu}:=\{ X\s A \mid |A|=\mu \}$.

\subsection{Organization of this paper}
In section 2, we introduce the notion of a \emph{stable poset},
observe that a stable poset witnesses the existence of an antichain,
and prove that under a very weak hypothesis, any poset of singular cofinality contains a stable poset,
concluding that this weak hypothesis implies the Milner-Sauer conjecture.

In section 3, we recall Hajnal and Sauer's definition of an \emph{antichain sequence},
and prove a somewhat surprising result: a poset of singular confinality contains a stable poset iff it contains an antichain sequence.
It is then an immediate corollary that our weak hypothesis indeed implies that
any poset of cofinality $\lambda>\cf(\lambda)=\kappa$ contains $\lambda^\kappa$ antichains of size $\kappa$.

\section{Existence of an antichain}
\begin{lemma}\label{14} Suppose $\langle P,\le\rangle$ is a poset, and $A\s B\s P$, then:
\begin{itemize}
\item[(a)] $\cf_P(\underline A)=\cf_P(A)\le |A|$.
\item[(b)] $\cf_P(A)\le\cf_P(B)\le\cf_P(\overline{B})$.
\item[(c)] $\cf_P\big(\bigcup_{\alpha<\mu}A_\alpha\big)\le\sum_{\alpha<\mu}\cf_P(A_\alpha)$ for any family $\{ A_\alpha\s P \mid \alpha<\mu\}$.
\item[(d)] If $\cf_P(B)>\cf_P(A)$, then $\cf_P(B\bks A)=\cf_P(B)$.
\end{itemize}
\end{lemma}
\begin{proof} Left as a warm-up exercise to the reader.\end{proof}

\begin{defn} Assume $\langle P,\leq \rangle$ is a poset of cofinality $\lambda>\cf(\lambda)=\kappa$.

A subset $P'\in[P]^\lambda$ is said to be \emph{stable} iff $\cf_P(P'\bks\overline{X})=\lambda$ for all $X\in[P']^{<\kappa}$.
\end{defn}

\begin{lemma}\label{23} Assume $\langle P,\leq \rangle$     is a poset of cofinality $\lambda>\cf(\lambda)=\kappa$.

If $P$ has a stable subset, then $P$ contains an antichain of size $\kappa$.
\end{lemma}
\begin{proof} Fix a stable subset $P'\s P$. We build an antichain $\{ x_\alpha\mid \alpha<\kappa\}\s P'$ by induction on $\alpha<\kappa$.
Suppose $X:=\{ x_\beta\mid\beta<\alpha\}\s P'$ have already been defined.
Since $X\in [P']^{<\kappa}$, $\cf_P(P'\bks\overline{X})=\lambda$.
Since $\cf_P(\underline{X})\le|X|<\kappa$,
we may find $x_\alpha\in P'$ such that $x_\alpha\not\in(\overline{X}\cup\underline{X})$.
\end{proof}

The notion of a stable subset gives rise to the following definition:
\begin{defn}\label{25} For a poset $\langle P,\le\rangle$ and a cardinal $\mu$, let:
$$\pd(P,\mu):=\{ X\in[P]^{<\mu}\mid \cf_P(P\bks\overline{X})< \cf_P(P)\}.$$
\end{defn}

Notice that if $X\in\pd(P,\mu)$ and $A\supseteq X$, then $A\in\pd(P,\theta)$ for all $\theta>|A|$.

\begin{lemma}\label{fact6} Assume $\langle P,\le\rangle$ is a poset of cofinality $\lambda>\cf(\lambda)=\kappa$.

For any $\mathcal A\s\pd(P,\lambda)$ of cardinality $\le\lambda$, there exists $Y\in[P]^{\kappa}$
such that $\underline{Y}\cap A\not=\emptyset$ for all $A\in\mathcal A$.
\end{lemma}
\begin{proof}
Let $\mathcal A=\{ X_i \mid i<\lambda\}$ be like in the hypothesis.
We shall construct a kind of anti $\mathcal A$-Luzin set.
Fix a strictly increasing sequence of cardinals converging to $\lambda$, $\langle \lambda_\alpha\mid \alpha<\kappa\rangle$.

Let $\alpha<\kappa$.
Put $\mathcal B_\alpha:=\{ X_i\mid i<\lambda_\alpha, \cf_P(P \bks\overline{X_i})<(\lambda_\alpha)^+\}$.
By $|\mathcal B_\alpha|<(\lambda_\alpha)^+$ and regularity of the latter, we have
 $\cf_P(P\bks\bigcap_{X\in\mathcal B_\alpha}\overline{X})=\cf_P(\bigcup_{X\in\mathcal B_\alpha}(P\bks\overline{X}))<(\lambda_\alpha)^+$.
Since $\cf_P(P)>(\lambda_\alpha)^+$, we may pick $y_\alpha\in \bigcap_{X\in\mathcal B_\alpha}\overline{X}$.

Finally, let $Y:=\{ y_\alpha\mid \alpha<\kappa\}$.
Since $\mathcal A=\bigcup_{\alpha<\kappa}\mathcal B_\alpha$, for each $X\in\mathcal A$,
there exists some $\alpha<\kappa$ with $y_\alpha\in Y\cap\overline{X}$, and hence, $\underline{Y}\cap X\not=\emptyset$.
\end{proof}

\begin{thm}\label{31} Assume $\langle P,\le\rangle$ is a poset of cofinality $\lambda>\cf(\lambda)=\kappa$.

For any cardinal $\mu<\lambda$ and $\mathcal A\s[P]^{<\mu}$ with $|\mathcal A|\le\lambda$,
there exists a subset $P'\s P$ such that:
\begin{itemize}
\item[(a)] $P'=\overline{P'}$ and $\cf_P(P\bks P')\le\mu\cdot\kappa$;
\item[(b)] $\cf_P(P'\bks \overline{A\cap P'})=\lambda$ for all $A\in\mathcal A$.
\end{itemize}
\end{thm}
\begin{proof}
For a set $Z\s P$, denote $\mathcal A_Z:=\pd(P,\mu)\cap\{ A\bks  Z\mid A\in\mathcal A\}$.
We define by induction on $\alpha<\mu$, a sequence of 
sets
$\langle Y_\alpha\in[P]^{\le\kappa}\mid \alpha<\mu\rangle$.

Assume $\langle Y_\beta\mid \beta<\alpha\rangle$ has already been defined.
Let $Z_\alpha:=\bigcup_{\beta<\alpha}\underline{Y_\beta}$ (where $Z_0:=\emptyset$).
By applying to Lemma \ref{fact6}, we pick $Y_\alpha\in[P]^{\kappa}$
such that $\underline{Y_\alpha}\cap A'\not=\emptyset$ for all $A'\in\mathcal A_{Z_\alpha}$.
This ends the construction.
Let $Y:=\bigcup_{\alpha<\mu} \underline{Y_\alpha}$ and $P':=P\bks Y$.
Evidently, $P'=\overline{P'}$ and $\cf_P(Y)\le\sum_{\alpha<\mu}|Y_\alpha|\le\mu\cdot\kappa$.

Towards a contradiction, assume there is $A\in\mathcal A$ with $\cf_P(P'\bks\overline{A\cap P'})<\lambda$.
It follows that $\cf_P(P\bks\overline{A\cap P'})\le\cf_P(Y)+\cf_P(P'\bks\overline{A\cap P'})<\lambda$,
that is, $(A\cap P')\in\pd(P,\mu)$.

We now define a function $f:\mu\rightarrow A$.
Fix $\alpha<\mu$.
Since $Z_\alpha\s Y$,  we have $A\cap P'\s A\bks Z_\alpha$ and it
follows from the remark after Definition \ref{25}, that $A\bks Z_\alpha\in\mathcal A_{Z_\alpha}$.
Thus, by the choice of $Y_\alpha$, we may pick some $f(\alpha)\in\underline{Y_\alpha}\cap(A\bks Z_\alpha)$.
This completes the definition of $f$.

Clearly, $f$ is an injection, and in particular, $|A|\ge\mu$,
contradicting the fact that $A\in\mathcal A\s[P]^{<\mu}$.
\end{proof}

\begin{defn}[Shelah \cite{Sh400}]\label{covdef}
For cardinals $\lambda\ge\kappa\ge\sigma>1, \mu\ge\kappa+\aleph_0$,
let: \[\cov(\lambda,\mu,\kappa,\sigma):=\min\{ |D| \mid D\s [\lambda]^{<\mu},\forall A\in [\lambda]^{<\kappa}\exists B\in [D]^{<\sigma}(A\s\bigcup B) \}. \]
\end{defn}
Thus, if $\lambda>\kappa$ are cardinals, then $\cf([\lambda]^{<\kappa},\s)=\cov(\lambda,\kappa,\kappa,2)$.
\begin{cor}\label{32} Assume $\cov(\lambda,\mu,\kappa,2)=\lambda$ for cardinals $\lambda>\mu\ge\cf(\lambda)=\kappa$.

If $\langle P,\le\rangle$ is a poset of cofinality $\lambda$, then $P$
contains a stable subset.
\end{cor}
\begin{proof} By restricting to a cofinal subset, we may assume that $|P|=\lambda$.
By   $\cov(\lambda,\mu,\kappa,2)=\lambda$, let us take
$\mathcal A\s[P]^{<\mu}$ such that $|\mathcal A|=\lambda$,
 and for each $X\in [P]^{<\kappa}$, there is $A\in\mathcal A$ with $X\s A$.
By applying to Theorem \ref{31}, we find $P'\s P$ such that $\cf_P(P'\bks \overline{A\cap P'})=\lambda$ for all $A\in\mathcal A$.
In particular, $\cf_P(P'\bks \overline{A'})=\lambda$ for all $A'\in [P']^{<\kappa}$, concluding that $P'$ is a stable subset.
\end{proof}

\begin{cor}\label{33} Assume $\cov(\lambda,\mu,\kappa,2)=\lambda$ for cardinals $\lambda>\mu\ge\cf(\lambda)=\kappa$.

If $\langle P,\le\rangle$ is a poset of cofinality $\lambda$, then $P$
contains an antichain of size~$\kappa$.
\end{cor}
\begin{proof} By Corollary \ref{32} and Lemma \ref{23}.
\end{proof}

\section{Antichain Sequences}
\begin{defn}[Hajnal-Sauer \cite{HaSa}]\label{331} Assume $\langle P,\leq \rangle$ is a poset,
and $\mathcal{A}=\langle A_{\alpha} \mid \alpha < \kappa \rangle$ is a family of mutually disjoint subsets of $P$.

$\mathcal{A}$ is said to be an \emph{antichain sequence} iff:
\begin{itemize}
\item[(a)] For all $\beta<\alpha<\kappa$, $|A_\beta|\le |A_\alpha|$;
\item[(b)] Any  $X\s\bigcup_{\alpha<\kappa}A_\alpha$ such that $|X\cap A_\alpha|\le 1$ for all $\alpha<\kappa$, is an antichain.
\end{itemize}
$\kappa$ is considered to be the \emph{length} of the antichain sequence,
and $\cf_P(\bigcup_{\alpha<\kappa}A_\alpha)$ as the \emph{cofinality} of the antichain sequence $\mathcal{A}$.
\end{defn}
It is worth noting that (b) is equivalent to the following statement:

(b*) For all $\beta<\alpha<\kappa$, $A_\alpha\cap\underline{A_\beta}=A_\alpha\cap\overline{A_\beta}=\emptyset$.

\begin{obs}\label{fact4} If $\langle P,\le\rangle$ is a poset of cofinality $\lambda>\cf(\lambda)=\kappa$,
and $P$ has an antichain sequence of length $\kappa$ and cofinality $\lambda$,
then $P$ contains $\lambda^\kappa$ antichains of size $\kappa$.
\end{obs}
\begin{proof} Fix $\mathcal{A}=\langle A_\alpha\mid \alpha<\kappa \rangle$ like in the hypothesis. For all $\alpha<\kappa$, set $\lambda_\alpha=|A_\alpha|$.
Finally, since $\langle \lambda_\alpha \mid \alpha<\kappa \rangle$ is non-decreasing, cofinal in $\lambda$:
$$|\{ \im(f) \mid f\in \prod_{\alpha<\kappa}A_\alpha\}|=\prod_{\alpha<\kappa}\lambda_\alpha=\lambda^\kappa.\qedhere$$
\end{proof}
We now aim at showing that the existence of an antichain sequence is equivalent to the existence of a stable subset.
To prove this, we first need the following essential observation.
\begin{lemma}[Hajnal-Sauer \cite{HaSa}]\label{spec_club} Assume $\langle P,\leq \rangle$ is a poset, and $P'\s P$.

If  $\cf_P(P')=|P'|=\lambda>\cf(\lambda)$, then $\sup\{ \cf_P(A)\mid A\in[P']^{<\lambda}\}=\lambda$.
\end{lemma}
\begin{proof} 
Put $\kappa:=\cf(\lambda)$.
By $|P'|=\lambda$, there exists a family of subsets $\{ A_\alpha\in[P']^{<\lambda}\mid\alpha<\kappa\}$
such that $P'=\bigcup_{\alpha<\kappa}A_\alpha$.
Let $\mu:=\sup\{ \cf_P(A_\alpha)\mid \alpha<\kappa\}$. If $\mu<\lambda$, then we obtain the following contradiction:
$$\lambda=\cf_P(P')=\cf_P\big(\bigcup_{\alpha<\kappa}A_\alpha\big)\le\sum_{\alpha<\kappa}\cf_P(A_\alpha)\le\kappa\cdot\mu<\lambda.\qedhere$$
\end{proof}

The main ideas of the following proof may already be found in \cite{rinot2}.
However, the following proof is simpler and more direct,
mainly due to new notion of a stable poset, who did not appear in \cite{rinot2}.
\begin{thm}\label{antiupwards} Assume $\langle P,\leq \rangle$ is a poset of cofinality $\lambda>\cf(\lambda)=\kappa$.

The following are equivalent:
\begin{itemize}
\item[(a)] There exists $\mathcal A=\langle A_\alpha\mid \alpha<\kappa\rangle$
 with $\{ A_\alpha\mid \alpha<\kappa\}\s[P]^{<\lambda}$
 such that $\mathcal A$ is an antichain sequence of length $\kappa$ and cofinality $\lambda$;
\item[(b)] There exists a stable subset $P'\s P$;
\item[(c)] There exists $Q\in [P]^\lambda$ with $\cf_P(Q)=\lambda$ such that $\cf_P(\overline{\{x \}}\cap Q)<\lambda$ for all $x\in Q$.
\end{itemize}
\end{thm}
\begin{proof}
Let $\langle \lambda_\alpha \mid \alpha <\kappa \rangle $ be a strictly-increasing sequence of cardinals converging to $\lambda$.
We prove (a)$\gorer$(b)$\gorer$(c)$\gorer$(a).

(a)$\gorer$(b) Suppose $\mathcal A$
is like in the hypothesis.
Put $P':=\bigcup_{\alpha<\kappa}A_\alpha$. By hypothesis, $\cf_P(P')=\lambda$,
and in particular $|P'|=\lambda$. Fix $X\in[P']^{<\kappa}$. By regularity of $\kappa$,
there exists some $\gamma<\kappa$ such that $X\s\bigcup_{\beta<\gamma}A_\beta$.
Since $\mathcal A$ is an antichain sequence, we get that $A_\delta\bks\overline{X}=A_\delta$ whenever $\gamma<\delta<\kappa$.
Since $|\bigcup_{\beta\le\gamma}A_\beta|<\lambda$, we must conclude that $\cf_P(P'\bks\overline{X})=\lambda$,
and hence $P'$ is a stable subset of $P$.

(b)$\gorer$(c) Let $P'\s P$ be stable and
assume towards a contradiction that:
$$(\star)\ (Q\in [P]^\lambda\wedge\cf_P(Q)=\lambda)\Rightarrow \exists x\in Q(\cf_P(\overline{\{x \}}\cap Q)=\lambda).$$
We build the following objects by induction on $\alpha<\kappa$:
\begin{itemize}
\item[(i)]  A set $\{ x_\alpha \mid \alpha<\kappa \}\s P'$.
\item[(ii)] A family of sets of the form $\{ A_{\alpha} \in [\overline{\{x_\alpha\}}\cap P']^{<\lambda} \mid \alpha < \kappa\}$.
\end{itemize}

Induction base: By $\cf_P(P')=\lambda$ and property $(\star)$,
we may pick $x_0\in P'$ such that $\cf_{P}(\overline{\{ x_0 \}}\cap P')=\lambda$, hence,
by Lemma \ref{spec_club} there exists $A_0\in[\overline{\{ x_0 \}}\cap P']^{<\lambda}$
with $\cf_{P}(A_0)>\lambda_0$.

Inductive step: Assume $X_\alpha:=\{ x_\beta\mid \beta<\alpha\}$ and $\{ A_\beta\mid \beta<\alpha\}$ have already been defined.
Since $P'$ is stable and $X_\alpha\in[P']^{<\kappa}$, we have that $\cf_P(P'\bks\overline{X_\alpha})=\lambda$.
It follows from $(\star)$, that we may choose $x_\alpha\in (P'\bks\overline{X_\alpha})$
such that $\cf_P(\overline{\{x_\alpha\}}\cap(P'\bks\overline{X_\alpha}))=\lambda$.
Thus, by applying to Lemma \ref{spec_club}, we pick
 $A_\alpha\in [\overline{\{x_\alpha\}}\cap(P'\bks\overline{X_\alpha})]^{<\lambda}$ with $\cf_P(A_\alpha)>\lambda_\alpha$.
 End of the construction.

Let $Q:=\bigcup_{\alpha<\kappa}A_\alpha$. Clearly, $\cf_P(Q)=\lambda$. Let $x\in Q$ be arbitrary.

To see that $\cf_P(\overline{\{x\}}\cap Q)<\lambda$, find $\alpha<\kappa$ with $x\in A_\alpha$.
In particular, $\overline{\{x\}}\s\overline{\{x_\alpha\}}\s\overline{X_{\alpha+1}}$,
and hence $\overline{\{x\}}\cap A_\delta
=\emptyset$ whenever $\alpha<\delta<\kappa$.
It follows that $(\overline{\{x\}}\cap Q)\s \bigcup_{\beta\le\alpha}A_\beta$ and $\cf_P(\overline{\{x\}}\cap Q)<\lambda$.
A contradiction.

(c)$\gorer$(a) Let $Q=\{ x_i \mid i<\lambda \}$ be like in the hypothesis.

Fix $\alpha<\kappa$ and set $B_\alpha:=\{ x_i\in Q \mid i<\lambda_\alpha, \cf_P(\overline{\{x_i\}}\cap Q)<\lambda_\alpha \}$. Thus:
$$\cf_P(\overline{B_\alpha}\cap Q)=\cf_P(\bigcup_{x\in B_\alpha}\overline{\{x\}}\cap Q)\le\sum_{x\in B_\alpha}\cf_P(\overline{\{x\}}\cap Q)\le\lambda_\alpha\cdot\lambda_\alpha=\lambda_\alpha.$$
Since $\{ B_\alpha \mid \alpha<\kappa \}$ is an increasing chain of sets, each of cardinality $<\lambda$,
and $\cf_P(\bigcup_{\alpha<\kappa}B_\alpha)=\cf_P(Q)=\lambda$, 
we may recursively define a strictly-increasing function $f:\kappa\rightarrow \kappa$, letting
$f(0):=\min\{\gamma <\kappa \mid \lambda_0<\cf_P(B_\gamma) \}$ and $f(\alpha):=\min\{\gamma <\kappa \mid \sum_{\beta<\alpha}\lambda_{f(\beta)}<\cf_P(B_\gamma) \}$
whenever $0<\alpha<\kappa$.

For all $\alpha<\kappa$, set $U_\alpha:=\bigcup_{\beta<\alpha}B_{f(\beta)}$ and $A_\alpha:=B_{f(\alpha)}\bks(\underline{U_\alpha}\cup\overline{U_\alpha})$.
To see that $\mathcal{A}:=\langle A_\alpha \mid \alpha<\kappa \rangle$ is an antichain sequence of cofinality $\lambda$,
we are left with showing that $\sup\{ \cf_P(A_\alpha) \mid \alpha<\kappa\}=\sup_{\alpha<\kappa}\lambda_\alpha$.

Fix $\alpha<\kappa$.
By $\cf_P(\underline{U_\alpha})\le \cf_P(\overline{U_\alpha}\cap Q)=\cf_P(\bigcup_{\beta<\alpha}\overline{B_{f(\beta)}}\cap Q)\le\sum_{\beta<\alpha}\lambda_{f(\beta)}$
and by the definition of $f$, we conclude that
$\cf_P(B_{f(\alpha)})>\cf_P((\underline{U_\alpha}\cup\overline{U_\alpha})\cap Q)$,
and hence $\cf_P(A_\alpha)=\cf_P(B_{f(\alpha)})$.
\end{proof}

\begin{cor}\label{33} Assume $\cov(\lambda,\mu,\kappa,2)=\lambda$ for cardinals $\lambda>\mu\ge\cf(\lambda)=\kappa$.

If $\langle P,\le\rangle$ is a poset of cofinality $\lambda$, then $P$
contains an antichain sequence of length $\kappa$ and cofinality $\lambda$.

In particular, every poset of cofinality $\lambda$ contains $\lambda^\kappa$ antichains of size $\kappa$.
\end{cor}
\begin{proof} By Corollary \ref{32}, Theorem \ref{antiupwards} and Observation \ref{fact4}.
\end{proof}

\section{Acknowledgements}
The author is grateful to M. Gitik for ``uncountably'' many discussions on the topic of this paper,
and for his valuable comments on this paper.

The author would also like to thank M. Pouzet and I. Gorelic for introducing him to the subject of
partial orders with singular cofinality.

\bibliographystyle{plain}

\end{document}